\documentclass[10pt]{article}
\usepackage{amssymb,latexsym}
\usepackage{epsfig}
\usepackage{eufrak}
\usepackage{amsmath}
\usepackage{mathrsfs}
\usepackage{color}

\newtheorem{theorem}{Theorem}
\newtheorem{lemma}{Lemma}

\newcommand{\be}{\begin{equation}}
\newcommand{\ee}{\end{equation}}
\newcommand{\bee}{\begin{eqnarray*}}
\newcommand{\eee}{\end{eqnarray*}}
\newcommand{\bel}{\begin{eqnarray}}
\newcommand{\eel}{\end{eqnarray}}
\newcommand{\bec}{\begin{cases}}
\newcommand{\eec}{\end{cases}}
\newcommand{\bem}{\begin{bmatrix}}
\newcommand{\eem}{\end{bmatrix}}

\newcommand{\la}{\label}
\newcommand{\li}{\left}
\newcommand{\ri}{\right}

\newcommand{\lc}{\lceil}
\newcommand{\rc}{\rceil}
\newcommand{\lf}{\lfloor}
\newcommand{\rf}{\rfloor}

\newcommand{\ga}{\gamma}
\newcommand{\Ga}{\Gamma}

\newcommand{\ka}{\kappa}

\newcommand{\f}{\frac}

\newcommand{\cd}{\cdots}

\newcommand{\qu}{\quad}
\newcommand{\qqu}{\qquad}

\newcommand{\mscr}{\mathscr}
\newcommand{\mcal}{\mathcal}

\newcommand{\bb}{\mathbb}

\newcommand{\mrm}{\mathrm}
\newcommand{\bs}{\boldsymbol}

\newcommand{\tx}{\text}

\newcommand{\iy}{\infty}

\newcommand{\bed}{\begin{description}}
\newcommand{\eed}{\end{description}}
\newcommand{\bei}{\begin{itemize}}
\newcommand{\eei}{\end{itemize}}
\newcommand{\ben}{\begin{enumerate}}
\newcommand{\een}{\end{enumerate}}
\newcommand{\bib}{\bibitem}
\newcommand{\beL}{\begin{lemma}}
\newcommand{\eeL}{\end{lemma}}
\newcommand{\beT}{\begin{theorem}}
\newcommand{\eeT}{\end{theorem}}
\newcommand{\sect}{\section}

\newcommand{\bpf}{\begin{pf}}
\newcommand{\epf}{\end{pf}}
\newcommand{\bsk}{\bigskip}
\newcommand{\bi}{\binom}

\newcommand{\pfbox}{\hfill\mbox{$\Box$}}

\newenvironment{pf}{\paragraph*{Proof{\rm.}}}{\pfbox\bigskip}

\begin{document}

\title{{\bf An Urn Model Approach for Deriving Multivariate Generalized Hypergeometric Distributions}
\thanks{The author had been previously working with Louisiana
State University at Baton Rouge, LA 70803, USA, and is now with Department of Electrical Engineering, Southern University and A\&M College,
Baton Rouge, LA 70813, USA; Email: chenxinjia@gmail.com}}

\author{Xinjia Chen}

\date{September  2013}

\maketitle

\begin{abstract}

We propose new generalized multivariate hypergeometric distributions, which extremely resemble the classical multivariate hypergeometric
distributions.  The proposed  distributions are derived based on an urn model approach.  In contrast to existing methods, this approach does not
involve hypergeometric series.

\end{abstract}

\section{Introduction}

It is well known that the parameters of the classical multivariate hypergeometric distributions are restricted to be positive integers. To
eliminate such restrictions, some researchers had proposed various generalized multivariate hypergeometric distributions. Among the works in
this direction, we recall the extension of hypergeometric distribution proposed by Kemp and Kemp \cite{Kemp}. A comprehensive survey of the
relevant results is given in Johnson et. al. \cite{Johnsona, Johnsonb}.  The main idea of existing methods for generalizing the classical
multivariate hypergeometric distributions is to use the hypergeometric series.  The generalized distributions derived from these methods
typically lack simplicity and do not perfectly resemble the classical multivariate hypergeometric distributions. In view of this situation, we
would like to propose a new generalization of multivariate hypergeometric distributions.  The proposed generalized distributions have  simple
structures and extremely resemble the classical multivariate hypergeometric distributions.  In sharp contrast to existing methods, we completely
abandon the use of hypergeometric series.  We derive the new generalization based on an urn model. We demonstrate that by appropriately
constructing multivariate P\'{o}lya-Eggenberger distributions with an urn model,  the proposed generalized distributions can be obtained by
letting certain parameters of the P\'{o}lya-Eggenberger distributions to be infinity. On the other hand, the proposed generalized distributions
include the P\'{o}lya-Eggenberger distributions as special cases.

The remainder of the paper is organized as follows. In Section 2, we propose a multivariate generalized hypergeometric distribution. In Section
3, we propose a multivariate generalized inverse hypergeometric distribution. In Section 4, we justify the proposed multivariate generalized
hypergeometric distribution. In Section 5, we justify the proposed multivariate generalized inverse hypergeometric distribution. For easy
reference, we include the derivation of the well-known multivariate P\'{o}lya-Eggenberger distributions in Appendices A na B.

Throughout this paper, we shall use the following notations.  Let $\bb{R}$ denote the set of all real numbers.  Let $\bb{Z}$ denote the set of
all integers.  Let $\bb{Z}^+$ denote the set of all nonnegative integers. Let $\bb{N}$ denote the set of all positive integers.  We use the
notation $\bi{t}{k}$ to denote a generalized combinatoric number in the sense that
\[
\bi{t}{k} = \f{\prod_{\ell = 1}^k (t -\ell + 1)}{k!} = \f{ \Ga(t + 1) }{\Ga(k + 1) \; \Ga(t - k + 1)},  \qqu  \bi{t}{0} = 1,
\]
where $t$ is a real number and $k$ is a non-negative integer.     The other notations will be made clear as we proceed.

\section{Multivariate Generalized Hypergeometric Distribution}

 In probability theory, random variables $X_0, X_1, \cd, X_\ka$, where $\ka \in \bb{N}$,  are said to possess a multivariate hypergeometric
distribution if  \[
\Pr \{ X_i = x_i, \; i = 0, 1, \cd, \ka \}  = \bec 0 & \tx{if $x_i \geq 1 + C_i$ for some $i \in \{ 0, 1, \cd, \ka \}$},\\
\f{ \prod_{i=0}^k \bi{  C_i }{x_i}   }{  \bi{N} {n} } & \tx{otherwise} \eec
\]
where $C_i \in \bb{N}, \; x_i \in \bb{Z}^+$ for $i = 0, 1, \cd, \ka$ and
\[
N = \sum_{i = 0}^\ka C_i  \geq \sum_{i = 0}^\ka x_i = n.
\]
Actually, under mild restrictions, the multivariate hypergeometric distribution can be generalized by allowing $N$ and $C_i$ to be real numbers.
More formally, we say that random variables $X_0, X_1, \cd, X_n$ possess a {\it multivariate generalized  hypergeometric distribution} if \be
\la{GINHY}
\Pr \{ X_i = x_i, \; i = 0, 1, \cd, \ka \} = \bec 0 & \tx{if $x_i \geq 1 + \mcal{C}_i > 1$ for some $i \in \{ 0, 1, \cd, k \}$},\\
\f{ \prod_{i=0}^k \bi{  \mcal{C}_i }{x_i}   }{  \bi{\mcal{N}} {n} } & \tx{otherwise} \eec \ee where $\mcal{C}_i \in \bb{R}, \; x_i \in \bb{Z}^+$
for $i = 0, 1, \cd, \ka$ and
\[
\mcal{N} = \sum_{i = 0}^k \mcal{C}_i \neq 0,  \qqu \sum_{i = 0}^\ka x_i = n, \qqu \f{n - 1}{\mcal{N}} < 1, \qqu
 \f{\mcal{C}_i}{\mcal{N}} > 0 \qu \tx{for} \;  i = 0, 1, \cd, k.
\]
The means and variances of $X_i$ are
given, respectively, as \bel &  & \bb{E} [ X_i ] =  \f{ n \mcal{C}_i }{\mcal{N}},   \la{exp889a}\\
&   & \mrm{Var} (X_i) =  \f{ n \mcal{C}_i (\mcal{N} - \mcal{C}_i) (\mcal{N} - n) }{ \mcal{N}^2 (\mcal{N} - 1)  },  \la{exp889b} \eel for $i = 0,
1, \cd, k$.

The justification of the proposed multivariate generalized  hypergeometric distribution and  equations (\ref{exp889a}), (\ref{exp889b}) is given
in Section \ref{CDFMGHapp}.  The proposed distribution includes many important distributions as special cases. Clearly, the multivariate
hypergeometric distribution is obtained from the multivariate generalized hypergeometric distribution by restricting $\mcal{C}_i, \; i = 0, 1,
\cd, k$ and $\mcal{N}$ as positive integers. The multivariate negative hypergeometric distribution is obtained from the multivariate generalized
hypergeometric distribution by taking $\mcal{C}_i, \; i = 0, 1, \cd, k$ and $\mcal{N}$ as negative integers.  The multinomial distribution is
obtained from the multivariate generalized hypergeometric distribution by letting $\mcal{N} \to \iy$ under the constraint that
$\f{\mcal{C}_i}{\mcal{N}}, \; i = 0, 1, \cd, k$ converge to positive numbers sum to $1$. The multivariate P\'{o}lya-Eggenberger distribution
\cite{Steyn} can be accommodated  as a special case of the multivariate generalized hypergeometric distribution.

\section{Multivariate Generalized Inverse Hypergeometric Distribution}

In probability theory, random variables $X_1, \cd, X_\ka$, where $\ka \in \bb{N}$,  are said to possess a {\it multivariate inverse
hypergeometric distribution} if  \[
\Pr \{ X_i = x_i, \; i = 1, \cd, \ka \} = \bec 0 & \tx{if $n \geq 1 + N$ or $x_i \geq 1 + C_i$ for some $i \in \{ 1, \cd, k \}$},\\
\f{\ga}{n} \f{ \prod_{i=0}^k \bi{ C_i }{x_i}   }{  \bi{N} {n} } & \tx{otherwise} \eec
\]
where $C_i \in \bb{N}, \; x_i \in \bb{Z}^+$ for $i = 1, \cd, k$ and
\[
C_0 \in \bb{N}, \qqu \ga \in \bb{N}, \qqu x_0 = \ga \leq C_0, \qqu n = \sum_{i=0}^k x_i, \qqu N = \sum_{i = 0}^k C_i. \]

Actually,  under mild restrictions, the multivariate inverse hypergeometric distribution can be generalized by allowing $N$ and $C_i$ to be real
numbers. More formally, we say that random variables $X_1, \cd, X_\ka$ possess a {\it multivariate generalized inverse hypergeometric
distribution} if \be \la{geninvhhyper89}
\Pr \{ X_i = x_i, \; i = 1, \cd, \ka \} = \bec 0 & \tx{if $\f{n - 1}{\mcal{N}} \geq 1$ or $x_i \geq 1 + \mcal{C}_i > 1$ for some $i \in \{ 1, \cd, k \}$},\\
\f{\ga}{n} \f{ \prod_{i=0}^k \bi{ \mcal{C}_i }{x_i}   }{  \bi{\mcal{N}} {n} } & \tx{otherwise} \eec \ee where $\mcal{C}_i \in \bb{R}, \; x_i \in
\bb{Z}^+$ for $i = 1, \cd, k$,
\[
\mcal{C}_0 \in \bb{R}, \qqu \ga \in \bb{N}, \qqu x_0 = \ga, \qqu n = \sum_{i = 0}^k x_i, \qqu \mcal{N} = \sum_{i = 0}^k \mcal{C}_i \neq 0, \qqu
\f{\mcal{C}_0} {\mcal{N}}
> \f{\ga - 1}{\mcal{N}},
\]
and $\f{\mcal{C}_i}{\mcal{N}} > 0$ for $i = 1, \cd, k$.  The means of $X_i$ are given as \be \la{expinv} \bb{E} [ X_i ] = \f{\ga
\mcal{C}_i}{\mcal{C}_0}, \qqu i = 1, \cd, k. \ee

The justification of the proposed distribution and  (\ref{expinv}) is given in Section \ref{CDFMGHINVapp}.

\sect{Derivation of  Multivariate Generalized Hypergeometric Distribution} \la{CDFMGHapp}

To justify that (\ref{GINHY}) indeed defines  a distribution, it suffices to consider two cases as follows.

Case (I):  All $\mcal{C}_i, \; i = 0, 1, \cd, k$ are integers and thus $\mcal{N}$ is an integer.  In this case, (\ref{GINHY}) defines a
classical multivariate hypergeometric distribution.

Case (II): There exists at least one index $j \in \{0,1, \cd, \ka \}$ such that $\mcal{C}_j$ is not an integer.  Without loss of generality,
assume that $\mcal{C}_0$ is not an integer throughout the remainder of this section.   To justify that (\ref{GINHY}) indeed defines  a
distribution in this case, we need some preliminary results.

\beL \la{lem888a}

\[
a + (\ell - 1) c  > 0, \qqu \ell = 1, \cd, n \qqu \tx{for $a > \f{n-1}{1 - \f{n-1}{\mcal{N}}}$}.
\]

\eeL

\bpf

For simplicity of notations, define \[ g(\ell) = a + (\ell - 1) c = a - (\ell - 1) \li \lc \f{a}{\mcal{N}} \ri
 \rc \qu \tx{for $\ell = 1, \cd, n$}.
\]
Clearly, $g(1) = a >0$. Recall the assumption that $\mcal{N}$ is a nonzero real number such that $\f{n - 1}{\mcal{N}} < 1$.
 Hence,  \bee  g(n) & = & a - (n - 1) \li \lc \f{a}{\mcal{N}} \ri \rc  >  a - (n - 1) \li ( \f{a}{\mcal{N}} + 1 \ri )\\
& = & \li ( 1 - \f{n - 1}{\mcal{N}} \ri ) a - (n - 1) > \li ( 1 - \f{n - 1}{\mcal{N}} \ri ) \f{n-1}{1 - \f{n-1}{\mcal{N}}} - (n - 1) =  0. \eee
Since $g(\ell)$ is a linear function of $\ell$, it must be true that
\[
g(\ell) > 0, \qqu \ell = 1, \cd, n.
\]
This proves the lemma.

\epf

We need to define some quantities. Let $a$ be a positive integer such that \be \la{asp88}
 a > \f{\mcal{N}}{\mcal{C}_i}, \qqu i = 1, \cd, k.
\ee For given $(\mcal{C}_0, \mcal{C}_1, \cd, \mcal{C}_k)$ and $\mcal{N}$, define \be \la{asp88bb} c = - \li \lc \f{a}{\mcal{N}} \ri \rc, \qqu
a_i = \li \lf \f{a \mcal{C}_i}{\mcal{N}} \ri \rf \qu \tx{for} \qu i = 1, \cd, k, \ee \be \la{asp88cc}
 a_0 = a - \sum_{i=1}^k a_i.
\ee Note that $a_0, a_1, \cd, a_\ka$ and $c$ are actually functions of $a$. We will use these functions as parameters to construct an urn model.
Based on the definition of these functions, we have the following results.

\beL \la{lem888b} Let $\mscr{S}$ denote the set of tuples $(x_0, x_1, \cd, x_k)$ of non-negative integers $x_0, x_1, \cd, x_k$ such that
$\sum_{i = 0}^k x_i = n$ and that there is no $i \in \{ 0, 1, \cd, k \}$ satisfying $x_i \geq 1 + \mcal{C}_i > 1$. Let $\mscr{S}_a^*$ denote the
set of tuples $(x_0, x_1, \cd, x_k)$ of non-negative integers $x_0, x_1, \cd, x_k$ such that $\sum_{i = 0}^k x_i = n$ and that $a_i + (x_i - 1)
c  > 0$  for $i \in \{ 0, 1, \cd, k \}$ such that $x_i > 0$. Then,  for large enough $a > 0$,
\[
\mscr{S} = \mscr{S}_a^*.
\]

\eeL

\bpf  If $\mcal{N} < 0$, then $\mcal{C}_i < 0$ for $i = 0, 1, \cd, k$.  It follows that $\mscr{S}$ is the set of tuples $(x_0, x_1, \cd, x_k)$
of non-negative integers $x_0, x_1, \cd, x_k$ such that $\sum_{i = 0}^k x_i = n$.  Moreover, as a consequence of (\ref{asp88}) and $\mcal{N} <
0$, we have $a_i + (x_i - 1) c  \geq a_i > 0$ for $i \in \{ 0, 1, \cd, k \}$ such that $x_i > 0$.  This implies that $\mscr{S}_a^*$ is the set
of tuples $(x_0, x_1, \cd, x_k)$ of non-negative integers $x_0, x_1, \cd, x_k$ such that $\sum_{i = 0}^k x_i = n$.  Therefore, in the case of
$\mcal{N} < 0$, we have shown that $\mscr{S} = \mscr{S}_a^*$ for large enough $a > 0$. It remains to show $\mscr{S} = \mscr{S}_a^*$ for large
enough $a > 0$ under the assumption that $\mcal{N} > 0$.  We proceed as follows.

First, we need to show that $\mscr{S} \subseteq \mscr{S}_a^*$ holds provided that $a > 0$ is sufficiently large. For this purpose, it suffices
to show that for any $(x_0, x_1, \cd, x_k) \in \mscr{S}$,
\[ a_i + (x_i - 1) c  > 0 \qqu \tx{holds for $i \in \{ 0, 1, \cd, k \}$ such that $x_i > 0$}
\]
provided that $a > 0$ is large enough.   Since $\mcal{N}$ is positive,  for any $(x_0, x_1, \cd, x_k) \in \mscr{S}$, it must be true that
$\f{x_i - 1}{\mcal{N}} < \f{\mcal{C}_i}{\mcal{N}}$ for $i = 0, 1, \cd, k$ such that $x_i
> 0$. Hence,
\[
\lim_{a \to \iy} \f{a_i}{a} = \f{\mcal{C}_i}{\mcal{N}} > \f{x_i - 1}{\mcal{N}} = \lim_{a \to \iy} \f{1}{a}(x_i - 1) \li \lc \f{a}{\mcal{N}} \ri
\rc, \qu \tx{for $i \in \{ 0, 1, \cd, k \}$ such that $x_i > 0$},
\]
which implies that
\[
\f{a_i}{a} > \f{1}{a}(x_i - 1) \li \lc \f{a}{\mcal{N}} \ri \rc, \qu \tx{for $i \in \{ 0, 1, \cd, k \}$ such that $x_i > 0$}
\]
for large enough $a > 0$. That is, if $a > 0$ is sufficiently large, then, \[ a_i > (x_i - 1) \li \lc \f{a}{\mcal{N}} \ri \rc = - (x_i - 1) c
\qu \qu \tx{for $i \in \{ 0, 1, \cd, k \}$ such that $x_i > 0$}.
\]
This establishes that $\mscr{S} \subseteq \mscr{S}_a^*$ holds provided that $a > 0$ is sufficiently large.

Next, we need to show that $\mscr{S} \supseteq \mscr{S}_a^*$ holds for large enough $a > 0$. Let $(x_0, x_1, \cd, x_k) \in \mscr{S}_a^*$.  Then,
$a_i + (x_i - 1) c  > 0$ for $i \in \{ 0, 1, \cd, k \}$ such that $x_i > 0$.  This means that \be \la{recall} a_i > (x_i - 1) \li \lc
\f{a}{\mcal{N}} \ri \rc \qu \tx{for $i \in \{ 0, 1, \cd, k \}$ such that $x_i > 0$}. \ee We need to show that $x_i < 1 + \mcal{C}_i$ for all $i
\in \{ 0, 1, \cd, k \}$. Clearly, $x_i < 1 + \mcal{C}_i$ for $x_i = 0$. It remains to show that if $a > 0$ is large enough, then $x_i < 1 +
\mcal{C}_i$ for $i \in \{ 0, 1, \cd, k \}$ such that $x_i > 0$. Making use of (\ref{recall}) and the assumption that $\mcal{N}
> 0$, we have \be \la{special8996}
 a_i > (x_i - 1) \f{a}{\mcal{N}}  \qu \tx{for $i \in \{ 0, 1, \cd, k \}$ such
that $x_i > 0$}. \ee By (\ref{special8996}) and the definition of $a_i$, we have
\[
 \f{a \mcal{C}_i}{\mcal{N}} \geq a_i > (x_i - 1) \f{a}{\mcal{N}}  \qu \tx{for $i \in \{ 1, \cd, k \}$ such that $x_i > 0$},
\]
which implies that $x_i < 1 + \mcal{C}_i$ for $i \in \{ 1, \cd, k \}$ such that $x_i > 0$.  Moreover, as a special case of (\ref{special8996}),
we have that
\[
a_0 > (x_0 - 1) \f{a}{\mcal{N}}  \qu \tx{if $x_0 > 0$}.
\]
It follows that
\[
\lim_{a \to \iy} \f{a_0}{a} = \f{\mcal{C}_0}{\mcal{N}} \geq (x_0 - 1) \f{1}{\mcal{N}}  \qu \tx{if $x_0 > 0$}.
\]
This implies that $\mcal{C}_0 \geq x_0 - 1$ if $x_0 > 0$.  For $x_0 \geq 1 + \mcal{C}_0$ to be valid, we must have $x_0 = 1 + \mcal{C}_0$, which
implies that $\mcal{C}_0$ is an integer. This contradicts to the assumption that $\mcal{C}_0$ is not an integer.  Therefore,  $x_0 < 1 +
\mcal{C}_0$ if $x_0 > 0$.  This proves that $\mscr{S} \supseteq \mscr{S}_a^*$ if $a > 0$ is sufficiently large.  Thus, we have shown that
$\mscr{S} = \mscr{S}_a^*$ for large enough $a > 0$ under the assumption that $\mcal{N} > 0$.  The proof of the lemma is thus completed.

\epf

To justify that (\ref{GINHY}) indeed defines  a distribution, it suffices to show that $\sum_{(x_0, x_1, \cd, x_k) \in \mscr{S} } P (x_1, \cd,
x_k) = 1$, where \[ P (x_1, \cd, x_k) = \f{ \prod_{i=0}^k \bi{  \mcal{C}_i }{x_i}   }{  \bi{\mcal{N}} {n} }.
\]  For this purpose, we use an urn model approach.  Assume that $a_0, a_1,
\cd, a_k$ and $c$ are functions of $a$ as defined by (\ref{asp88}),  (\ref{asp88bb}) and (\ref{asp88cc}).  Consider an urn containing $a_0, a_1,
\cd, a_k$ initial balls of $k+1$ different colors, $\bb{C}_0, \bb{C}_1, \cd, \bb{C}_k$, respectively. The sampling scheme is as follows.  A ball
is drawn at random from the urn. The color of the drawn ball is noted and then the ball is returned to the urn along with $c$ additional balls
of the same color. In the case that after drawing a ball, the number of balls of that color remained in the urn is no greater than $- (c + 1)$,
that type of balls will be eliminated from the sampling experiment. This operation is repeated, using the newly constituted urn, until $n$ such
operations have been completed.     Let $X_0^*, X_1^*, X_2^*, \cd, X_k^*$ denotes the numbers  of balls of colors $\bb{C}_0, \bb{C}_1, \bb{C}_2,
\cd, \bb{C}_k$, respectively, drawn at the end of $n$ trials. From Lemma \ref{lem888a}, we have that $\prod_{\ell = 1}^n [a + (\ell - 1) c ]
> 0$ for large enough $a > 0$. Let $(x_0, x_1, \cd, x_k) \in \mscr{S}_a^*$ with $a > 0$ large enough.
It can be shown that $\Pr \{ X_1^* = x_1, X_2^* = x_2, \cd, X_k^* = x_k \} = P^* (x_1, \cd, x_k)$, where {\small \bee P^* (x_1, \cd, x_k) =
\f{1}{ \prod_{\ell = 1}^n [a + (\ell - 1) c ] } \bi{n} {\bs{x}} \prod_{i=0}^k \li \{ \prod_{\ell = 1}^{x_i} [a_i + (\ell - 1) c ] \ri \}, \eee}
where $\sum_{i=0}^k x_i = n$ and $\bi{n} {\bs{x}} = \f{n!}{x_0! x_1!  \cd x_k!}$ is the multinomial coefficient.  This is the well-known
P\'{o}lya-Eggenberger distribution (see, Appendix A for its derivation). Clearly,
\[
\sum_{(x_0, x_1, \cd, x_k) \in \mscr{S}_a^*} P^* (x_1, \cd, x_k) = 1
\]
for large enough $a > 0$.  Recalling Lemma \ref{lem888b},  we have  that for large enough $a > 0$,
\[
\mscr{S} = \mscr{S}_a^*,  \] where $\mscr{S}$ is independent of $a$.  This implies that
\[
\sum_{(x_0, x_1, \cd, x_k) \in \mscr{S} } P^* (x_1, \cd, x_k) = \sum_{(x_0, x_1, \cd, x_k) \in \mscr{S}_a^*} P^* (x_1, \cd, x_k) = 1
\]
for large enough $a > 0$.  Note that the number of all $(x_0, x_1, \cd, x_k)$-tuples in $\mscr{S}$ is finite and that \bee &  & P (x_1, \cd,
x_k) = \f{ 1 } { \prod_{\ell = 1}^n \li [ 1 - (\ell - 1) \f{1}{\mcal{N}} \ri ] } \bi{n} {\bs{x}} \prod_{i=0}^k
\li \{ \prod_{\ell = 1}^{x_i} \li [ \f{\mcal{C}_i}{\mcal{N}} - (\ell - 1) \f{1}{\mcal{N}} \ri ] \ri \},\\
&  & P^* (x_1, \cd, x_k)  =   \f{ 1 } { \prod_{\ell = 1}^n \li [ 1 + (\ell - 1) \nu \ri ] } \bi{n} {\bs{x}} \prod_{i=0}^k \li \{ \prod_{\ell =
1}^{x_i} \li [ p_i + (\ell - 1) \nu \ri ] \ri \},  \eee where $\nu = \f{c}{a}, \; p_0 = \f{ a_0 }{a}$ and $p_i = \f{a_i}{a} = \f{1}{a} \li \lf
\f{ a \mcal{C}_i}{\mcal{N}} \ri \rf$ for $i = 1, \cd, k$.  For any tuple in $\mscr{S}$, the probability $P^* (x_1, \cd, x_k)$ is a continuous
function of $\nu$ and $(p_0, p_1, \cd, p_k)$.  Since $\nu \to - \f{1}{\mcal{N}}$ and $p_i \to \f{\mcal{C}_i}{\mcal{N}}$ for $i = 0, 1, \cd, k$
as $a \to \iy$, it follows that $P^* (x_1, \cd, x_k) \to P (x_1, \cd, x_k)$ for any tuple $(x_0, x_1, \cd, x_k) \in \mscr{S}$ as $a \to \iy$.
Consequently, \[ \sum_{(x_0, x_1, \cd, x_k) \in \mscr{S} } P^* (x_1, \cd, x_k)  \to \sum_{(x_0, x_1, \cd, x_k) \in \mscr{S} } P (x_1, \cd, x_k)
\] as $a \to \iy$. Since $\sum_{(x_0, x_1, \cd, x_k) \in \mscr{S} } P^* (x_1, \cd, x_k)  = 1$ for large enough $a > 0$, it must be true that $\sum_{(x_0, x_1, \cd, x_k)
\in \mscr{S} } P (x_1, \cd, x_k) = 1$.   Thus, we have justified that (\ref{GINHY}) indeed defines  a distribution.

Note that the means and variances of $X_i^*$ are
given, respectively, by \bee &  & \bb{E} [ X_i^* ] =  \f{ n a_i }{a},   \\
&   & \mrm{Var} (X_i^*) =  \f{ n a_i (a - a_i) (a + n c) }{ a^2 (a + c)  }  \eee for $i = 0, 1, \cd, k$. Making use of these results and letting
$a \to \iy$ lead to (\ref{exp889a}) and (\ref{exp889b}).

\sect{Derivation of  Multivariate Generalized Inverse Hypergeometric Distribution} \la{CDFMGHINVapp}

 To justify that (\ref{geninvhhyper89}) indeed defines  a distribution, we need some preliminary results.

\beL \la{lem888ainv}

Let $x_0 = \ga$ and $x_1, \cd, x_k$ be nonnegative integers.  Let $n = \sum_{\ell = 0}^k x_\ell$.  Assume that $\f{n - 1}{\mcal{N}} < 1$. Then,
\[
 a + (\ell - 1) c  > 0, \qqu \ell = 1, \cd, n \qqu \tx{for $a > \f{n-1}{1 - \f{n-1}{\mcal{N}}}$}.
 \]

\eeL

\bpf

The lemma can be shown by using the same argument as that of Lemma \ref{lem888a} and the assumption that $\f{n - 1}{\mcal{N}} < 1$.

\epf

We need to define some quantities. Let $a > 0$ be a positive integer such that \be \la{asp88inv}
 a > \f{\mcal{N}}{\mcal{C}_i}, \qqu i = 1, \cd, k.
\ee For given $(\mcal{C}_0, \mcal{C}_1, \cd, \mcal{C}_k)$ and $\mcal{N}$, define \be \la{asp88invbb} c = - \li \lc \f{a}{\mcal{N}} \ri \rc, \qqu
a_i = \li \lf \f{a \mcal{C}_i}{\mcal{N}} \ri \rf \qu \tx{for} \qu i = 1, \cd, k, \ee \be \la{asp88invcc} a_0 = a - \sum_{i=1}^k a_i. \ee

\beL \la{lem888binv} Let $\mscr{S}$ denote the set of tuples $(x_0, x_1, \cd, x_k)$, where  $x_0 = \ga$ and $x_1, \cd, x_k$ are non-negative
integers, such that $\sum_{i = 0}^k x_i = n$ and that there is no $i \in \{ 1, \cd, k \}$ satisfying $\f{n - 1}{\mcal{N}} \geq 1$ or $x_i \geq 1
+ \mcal{C}_i > 1$. Let $\mscr{S}_a^*$, with $a > 0$,  denote the set of tuples $(x_0, x_1, \cd, x_k)$, where  $x_0 = \ga$ and $x_1, \cd, x_k$
are non-negative integers, such that $\sum_{i = 0}^k x_i = n$ and that $a + (n-1) c > 0, \; a_i + (x_i - 1) c  > 0$  for $i \in \{ 1, \cd, k \}$
such that $x_i
> 0$.    Then, for large enough $a > 0$,
\[
\mscr{S} = \mscr{S}_a^*.
\]

\eeL

\bpf  If $\mcal{N} < 0$, then $\mcal{C}_i < 0$ for $i = 0, 1, \cd, k$. Hence, $\mscr{S}$ is the set of tuples $(x_0, x_1, \cd, x_k)$, where $x_0
= \ga$ and $x_1, \cd, x_k$ are non-negative integers,   such that $\sum_{i = 0}^k x_i = n$. Moreover, as a consequence of (\ref{asp88inv}) and
$\mcal{N} < 0$, we have $a_i + (x_i - 1) c  \geq a_i > 0$ for $i \in \{ 1, \cd, k \}$ such that $x_i > 0$.  It follows that $\mscr{S}_a^*$ is
the set of tuples $(x_0, x_1, \cd, x_k)$ of non-negative integers $x_0, x_1, \cd, x_k$ such that $\sum_{i = 0}^k x_i = n$. Therefore, $\mscr{S}
= \mscr{S}_a^*$ holds for large enough $a > 0$ in the case of $\mcal{N} < 0$.  It remains to show $\mscr{S} = \mscr{S}_a^*$ for large enough $a
> 0$ under the assumption that $\mcal{N} > 0$.

First, we need to show that $\mscr{S} \subseteq \mscr{S}_a^*$ holds provided that $a > 0$ is sufficiently large. For this purpose, it suffices
to show that for any $(x_0, x_1, \cd, x_k) \in \mscr{S}$,
\[ a_i + (x_i - 1) c  > 0 \qqu \tx{holds for $i \in \{ 1, \cd, k \}$ such that $x_i > 0$},
\]
provided that $a > 0$ is large enough.  Note that for any $(x_0, x_1, \cd, x_k) \in \mscr{S}$, it must be true that $\f{n - 1}{\mcal{N}} < 1$
and $\f{x_i - 1}{\mcal{N}} < \f{\mcal{C}_i}{\mcal{N}}$ for $i \in \{ 1, \cd, k \}$ such that $x_i
> 0$.  Since $\f{n - 1}{\mcal{N}} < 1$, it follows from Lemma \ref{lem888ainv} that $a + (n-1) c > 0$ for large enough $a > 0$.
Since $\f{x_i - 1}{\mcal{N}} < \f{\mcal{C}_i}{\mcal{N}}$ for $i \in \{ 1, \cd, k \}$ such that $x_i
> 0$, we have
\[
\lim_{a \to \iy} \f{a_i}{a} = \f{\mcal{C}_i}{\mcal{N}} > \f{x_i - 1}{\mcal{N}} = \lim_{a \to \iy} \f{1}{a}(x_i - 1) \li \lc \f{a}{\mcal{N}} \ri
\rc, \qu \tx{for $i \in \{ 1, \cd, k \}$ such that $x_i > 0$},
\]
which implies that
\[
\f{a_i}{a} > \f{1}{a}(x_i - 1) \li \lc \f{a}{\mcal{N}} \ri \rc, \qu \tx{for $i \in \{ 1, \cd, k \}$ such that $x_i > 0$}
\]
for large enough $a > 0$. That is, if $a > 0$ is sufficiently large, then, \[ a_i > (x_i - 1) \li \lc \f{a}{\mcal{N}} \ri \rc \qu \qu \tx{for $i
\in \{ 1, \cd, k \}$ such that $x_i > 0$}.
\]
This establishes that $\mscr{S} \subseteq \mscr{S}_a^*$ holds provided that $a > 0$ is sufficiently large.

Next, we need to show that $\mscr{S} \supseteq \mscr{S}_a^*$ holds for large enough $a > 0$.  Let $(x_0, x_1, \cd, x_k) \in \mscr{S}_a^*$ for $a
> 0$. Then, $a + (n-1) c > 0$ and $a_i + (x_i - 1) c  > 0$ for $i \in \{ 1, \cd, k \}$ such that $x_i > 0$.  As a consequence of $a + (n-1) c >
0$, we have that for $a > 0$,
\[
 a - (n-1)  \f{a}{\mcal{N}}  \geq  a - (n-1) \li \lc \f{a}{\mcal{N}} \ri \rc = a + (n-1) c > 0,
\]
which implies that $\f{n - 1}{\mcal{N}} < 1$.  Since $a_i + (x_i - 1) c  > 0$ for $i \in \{ 1, \cd, k \}$ such that $x_i > 0$, it follows that
\be \la{recallinv} a_i > (x_i - 1) \li \lc \f{a}{\mcal{N}} \ri \rc \qu \tx{for $i \in \{ 1, \cd, k \}$ such that $x_i > 0$}. \ee We need to show
that $x_i < 1 + \mcal{C}_i$ for $i \in \{ 1, \cd, k \}$. Clearly, $x_i < 1 + \mcal{C}_i$ for $x_i = 0$. It remains to show that if $a > 0$ is
large enough, then $x_i < 1 + \mcal{C}_i$ for $i \in \{ 1, \cd, k \}$ such that $x_i > 0$. Making use of (\ref{recallinv}) and the assumption
that $\mcal{N} > 0$, we have \be \la{special8996inv}
 a_i > (x_i - 1) \f{a}{\mcal{N}}  \qu \tx{for $i \in \{ 1, \cd, k \}$ such
that $x_i > 0$}. \ee By virtue of (\ref{special8996inv}) and the definition of $a_i$, we have
\[
 \f{a \mcal{C}_i}{\mcal{N}} \geq a_i > (x_i - 1) \f{a}{\mcal{N}}  \qu \tx{for $i \in \{ 1, \cd, k \}$ such that $x_i > 0$},
\]
which implies that $x_i < 1 + \mcal{C}_i$ for $i \in \{ 1, \cd, k \}$ such that $x_i > 0$. This proves that $\mscr{S} \supseteq \mscr{S}_a^*$ if
$a > 0$ is sufficiently large.  Thus, we have shown that $\mscr{S} = \mscr{S}_a^*$ for large enough $a > 0$.

\epf

\beL \la{unboundeda}
\[
a_0 + (x_0 - 1) c > 0
\]
for large enough $a > 0$. \eeL

\bpf By the assumption that $\f{ \mcal{C}_0 }{ \mcal{N} } > \f{\ga - 1}{ \mcal{N} }$, we have that $\f{ a \mcal{C}_0 }{ \mcal{N} } > 1 + \f{(\ga
- 1) a}{ \mcal{N} }$ for large enough $a > 0$, which implies that
\[
a - \sum_{i=1}^k \f{ a \mcal{C}_i }{ \mcal{N} }  > 1 +  \f{(\ga - 1) a}{ \mcal{N} }
\]
for large enough $a > 0$.  Therefore,
\[
a - \sum_{i=1}^k \li \lf \f{ a \mcal{C}_i }{ \mcal{N} } \ri \rf  >  \li \lc \f{(\ga - 1) a}{ \mcal{N} } \ri \rc
\]
for large enough $a > 0$.  Recalling the definitions of $a_0, x_0$ and $c$, we have that $a_0 + (x_0 - 1) c > 0$ for large enough $a > 0$. \epf

\bsk

To justify that (\ref{geninvhhyper89}) indeed defines  a distribution, it suffices to show that $\sum_{(x_0, x_1, \cd, x_k) \in \mscr{S} } P
(x_1, \cd, x_k) = 1$, where \[ P (x_1, \cd, x_k) = \f{\ga}{n} \f{ \prod_{i=0}^k \bi{  \mcal{C}_i }{x_i}   }{  \bi{\mcal{N}} {n} }
\]
with $x_0 = \ga$ and $n = \sum_{i=0}^k x_i$.  For this purpose, we use an urn model approach.  Assume that $a_0, a_1, \cd, a_k$ and $c$ are
functions of $a > 0$ as defined by (\ref{asp88inv}), (\ref{asp88invbb}) and (\ref{asp88invcc}).   For simplicity of notations, define $\mscr{A}
= \{a > 0: a_0 + (x_0 - 1) c > 0 \}$.  By Lemma \ref{unboundeda} and the definition of $\mscr{A}$, the value of $a \in \mscr{A}$ can be
arbitrarily large. Consider an urn containing $a_0, a_1, \cd, a_k$ initial balls of $k+1$ different colors, $\bb{C}_0, \bb{C}_1, \cd, \bb{C}_k$,
respectively. The sampling scheme is as follows.  A ball is drawn at random from the urn. The color of the drawn ball is noted and then the ball
is returned to the urn along with $c$ additional balls of the same color. In the case that after drawing a ball, the number of balls of that
color remained in the urn is no greater than $- (c + 1)$, that type of balls will be eliminated from the sampling experiment. This operation is
repeated, using the newly constituted urn, until $\ga$ balls of color $\bb{C}_0$ have been chosen. From Lemma \ref{lem888ainv}, we have that
$\prod_{\ell = 1}^n [a + (\ell - 1) c ] > 0$ for large enough $a > 0$.  For $(x_0, x_1, \cd, x_k) \in \mscr{S}_a^*$ with $a \in \mscr{A}$, it
must be true that $\prod_{i=0}^k \li \{ \prod_{\ell = 1}^{x_i} [a_i + (\ell - 1) c ] \ri \} > 0$.  Hence, for large enough $a \in \mscr{A}$, the
numbers $X_1^*, X_2^*, \cd, X_k^*$ of balls of colors $\bb{C}_1, \bb{C}_2, \cd, \bb{C}_k$, respectively, drawn at the end of $n = \sum_{i=0}^k
x_i$ trials have the joint probability mass function {\small \bee P^* (x_1, \cd, x_k) = \Pr \{ X_1^* = x_1, X_2^* = x_2, \cd, X_k^* = x_k \} =
\f{1}{ \prod_{\ell = 1}^n [a + (\ell - 1) c ] } \f{\ga}{n} \bi{n} {\bs{x}} \prod_{i=0}^k \li \{ \prod_{\ell = 1}^{x_i} [a_i + (\ell - 1) c ] \ri
\}, \eee} where  $\bi{n} {\bs{x}} = \f{n!}{x_0! x_1!  \cd x_k!}$ is the multinomial coefficient.  This is the well-known inverse
P\'{o}lya-Eggenberger distribution. Clearly,
\[
\sum_{(x_0, x_1, \cd, x_k) \in \mscr{S}_a^*} P^* (x_1, \cd, x_k) = 1
\]
for large enough $a \in \mscr{A}$.  Recalling Lemma \ref{lem888binv},  we have  that for large enough $a \in \mscr{A}$,
\[
\mscr{S} = \mscr{S}_a^*,  \] where $\mscr{S}$ is independent of $a$.  This implies that
\[
\sum_{(x_0, x_1, \cd, x_k) \in \mscr{S} } P^* (x_1, \cd, x_k) = \sum_{(x_0, x_1, \cd, x_k) \in \mscr{S}_a^*} P^* (x_1, \cd, x_k) = 1
\]
for large enough $a \in \mscr{A}$.  Note that the number of all $(x_0, x_1, \cd, x_k)$-tuples in $\mscr{S}$ is finite and that \bee &  & P (x_1,
\cd, x_k) = \f{ 1 } { \prod_{\ell = 1}^n \li [ 1 - (\ell - 1) \f{1}{\mcal{N}} \ri ] } \f{\ga}{n} \bi{n} {\bs{x}} \prod_{i=0}^k
\li \{ \prod_{\ell = 1}^{x_i} \li [ \f{\mcal{C}_i}{\mcal{N}} - (\ell - 1) \f{1}{\mcal{N}} \ri ] \ri \},\\
&  & P^* (x_1, \cd, x_k)  =   \f{ 1 } { \prod_{\ell = 1}^n \li [ 1 + (\ell - 1) \nu \ri ] } \f{\ga}{n} \bi{n} {\bs{x}} \prod_{i=0}^k \li \{
\prod_{\ell = 1}^{x_i} \li [ p_i + (\ell - 1) \nu \ri ] \ri \},  \eee where $\nu = \f{c}{a}, \; p_0 = \f{ a_0 }{a}$ and $p_i = \f{a_i}{a} =
\f{1}{a} \li \lf \f{ a \mcal{C}_i}{\mcal{N}} \ri \rf$ for $i = 1, \cd, k$.  For any tuple in $\mscr{S}$, the probability $P^* (x_1, \cd, x_k)$
is a continuous function of $\nu$ and $(p_0, p_1, \cd, p_k)$.  Since $\nu \to - \f{1}{\mcal{N}}$ and $p_i \to \f{\mcal{C}_i}{\mcal{N}}$ for $i =
0, 1, \cd, k$ as $a \in \mscr{A}$ tends to $\iy$, it follows that $P^* (x_1, \cd, x_k) \to P (x_1, \cd, x_k)$ for any tuple $(x_0, x_1, \cd,
x_k) \in \mscr{S}$ as $a \in \mscr{A}$ tends to $\iy$. Consequently, \[ \sum_{(x_0, x_1, \cd, x_k) \in \mscr{S} } P^* (x_1, \cd, x_k)  \to
\sum_{(x_0, x_1, \cd, x_k) \in \mscr{S} } P (x_1, \cd, x_k)
\] as $a \in \mscr{A}$ tends to $\iy$. Since $\sum_{(x_0, x_1, \cd, x_k) \in \mscr{S} } P^* (x_1, \cd, x_k)  = 1$ for large enough $a \in \mscr{A}$, it must be true that
$\sum_{(x_0, x_1, \cd, x_k) \in \mscr{S} } P (x_1, \cd, x_k) = 1$.   Thus, we have justified that (\ref{geninvhhyper89}) indeed defines  a
distribution.

Note that
\[
\bb{E} [ X_i^* ] = \f{n a_i}{a}, \qqu i = 1, \cd, k.
\]
Making use of this result and letting $a \in \mscr{A}$ tend to infinity lead to (\ref{expinv}).

\appendix

\section{Multivariate P\'{o}lya-Eggenberger Distribution}

Consider an urn containing $a_0, a_1, \cd, a_k$ initial balls of $k+1$ different colors, $\bb{C}_0, \bb{C}_1, \cd, \bb{C}_k$, respectively. The
sampling scheme is as follows.  A ball is drawn at random from the urn. The color of the drawn ball is noted and then the ball is returned to
the urn along with $c$ additional balls of the same color.  In the case that after drawing a ball, the number of balls of that color remained in
the urn is no greater than $- (c + 1)$, that type of balls will be eliminated from the sampling experiment.  This operation is repeated, using
the newly constituted urn, until $n$ such operations (often called ``trials'') have been completed. Assume that $a + (n-1) c
> 0$.  Steyn \cite{Steyn} showed that the  numbers $X_1, X_2, \cd, X_k$ of balls of colors $\bb{C}_1, \bb{C}_2, \cd, \bb{C}_k$,
respectively, drawn at the end of $n$ trials have the joint probability mass function {\small \bee \Pr \{ X_1 = x_1, \cd, X_k = x_k \}  = \bec
\f{\bi{n} {\bs{x}} \prod_{i=0}^k \prod_{\ell = 1}^{x_i} [a_i + (\ell - 1) c ]}{ \prod_{\ell = 1}^n [a + (\ell - 1) c ] }   & \tx{if $a_i + (\ell
- 1) c > 0$ for all $i \in \{1, \cd, k \}$ such
that $x_i > 0$}, \\
0 & \tx{otherwise} \eec \eee} where
\[ \sum_{i=0}^k x_i = n, \qqu \sum_{i=0}^k a_i = a,
\]
and \[ \bi{n} {\bs{x}} = \f{n!}{x_0! x_1!  \cd x_k!}
\]
is the multinomial coefficient.  This is the well-known multivariate P\'{o}lya-Eggenberger distribution.  The derivation of this distribution is
as follows.

To make the sampling experiment well-defined,  assume that for $i = 0, 1, \cd, k$, there are infinitely many balls of color $\bb{C}_i$.  Let
these balls be labeled as $B_{i, 1}, B_{i, 2}, \cd$ so that $B_{i, j}, \; j = 1, \cd, a_i$ will be initially put in the urn and that balls
$B_{i, j}, \; j = a_i + 1, \; a_i + 2, \cd$ are used as additional balls, which are used in an order consistent with their indexes, that is, the
ball with index $j$ must have been added to the urn if the ball with index $j + 1$ is to be added to the urn. Let the set of balls be denoted by
$\{ B_{i, j} \}$.

Note that every sequence of $n$ drawings can be represented by a permutation like $P_1 P_2 \cd P_n$ of length $n$, where $P_\ell$ denotes the
ball drawing in the $\ell$-th trial.   Clearly, $P_\ell$ is picked from a subset of $\{ B_{i, j} \}$ under certain constraints.   Note that the
same $B_{i, j}$ can appear in different position of $P_1 P_2 \cd P_n$.   For $\ell = 1, \cd, n$, in the $\ell$-th trial, there are $a + (\ell -
1) c$ balls in the urn for equally likely random drawing.  It follows that there are $\prod_{\ell = 1}^n [ a + (\ell - 1) c ]$ possible
permutations, each of them are equally likely.  For the experiment to be feasible, $\prod_{\ell = 1}^n [ a + (\ell - 1) c ]$ must be positive.

Now consider the making of permutation like $P_1 P_2 \cd P_n$ which corresponds to $x_0, x_1, \cd, x_k$ balls of color $\bb{C}_0, \bb{C}_1, \cd,
\bb{C}_k$, respectively.  Such permutations can be made as follows.

First, for $i = 0, 1, \cd, k$, we choose $x_i$ numbers from $\{1, 2, \cd, n \}$ and denote the set of the numbers as $S_i$. Thus, we have $k$
sets $S_0, S_1, \cd, S_k$. There are $\f{n!}{x_0 ! x_1 !, \cd, x_k !}$ ways to create such sets.

Second, for $i = 0, 1, \cd, k$, arrange the numbers in $S_i$ as $r_{\ell, i}, \; \ell = 1, \cd, x_i$.  In the $r_{\ell, i}$-th trials, there are
$a_i + (\ell - 1) c$ balls of color $\bb{C}_i$.  Hence, for $i = 0, 1, \cd, k$, there are $\prod_{\ell = 1}^{x_i} [ a_i + (\ell - 1) c ]$ ways
to make a permutation having $x_i$ balls if the pattern of balls of other colors are fixed.  It follows that there are $\bi{n} {\bs{x}}
\prod_{i=0}^k \li \{ \prod_{\ell = 1}^{x_i} [a_i + (\ell - 1) c ] \ri \}$ permutations,  corresponding to $x_i$ balls of color $\bb{C}_i$ for $i
= 0, 1, \cd, k$.  For the experiment to be feasible, $\prod_{i=0}^k \li \{ \prod_{\ell = 1}^{x_i} [a_i + (\ell - 1) c ] \ri \}$ must be
positive.

Third, since each permutation is equally likely, the probability of getting $x_i$ balls of color $\bb{C}_i$ respectively is equal to  the ratio
of the number of permutations,  having $x_i$ balls of color $\bb{C}_i$ for $i = 0, 1, \cd, k$,  to the total number of permutations, that is,
\[
\f{1}{ \prod_{\ell = 1}^n [a + (\ell - 1) c ] } \bi{n} {\bs{x}} \prod_{i=0}^k \li \{ \prod_{\ell = 1}^{x_i} [a_i + (\ell - 1) c ] \ri \}.
\]
This completes the proof of the P\'{o}lya-Eggenberger distribution.

\section{Multivariate Inverse P\'{o}lya-Eggenberger Distribution}

Consider an urn containing $a_0, a_1, \cd, a_k$ initial balls of $k+1$ different colors, $\bb{C}_0, \bb{C}_1, \cd, \bb{C}_k$, respectively. The
sampling scheme is as follows.  A ball is drawn at random from the urn. The color of the drawn ball is noted and then the ball is returned to
the urn along with $c$ additional balls of the same color. In the case that after drawing a ball, the number of balls of that color remained in
the urn is no greater than $- (c + 1)$, that type of balls will be eliminated from the sampling experiment. This operation is repeated, using
the newly constituted urn, until $\ga$ balls of color $\bb{C}_0$ have been chosen.  Assume that $a_0 + c (\ga - 1) > 0$.  The joint distribution
of the numbers $X_1, \cd, X_k$ of balls of colors $\bb{C}_1, \cd, \bb{C}_k$, respectively, drawn when this requirement is achieved is  {\small
\[ \Pr \{ X_1 = x_1, \cd, X_k = x_k  \} = \bec \f{\ga \bi{n}{\bs{x}} \prod_{i=0}^k \prod_{\ell = 1}^{x_i} [a_i + (\ell - 1) c ] }{ n \prod_{\ell
= 1}^{n} [a + (\ell - 1) c ] }  & \tx{  if $a_i + (x_i - 1) c >
0$ for $i \in \{1, \cd, k \}$ such that $x_i > 0$  }, \\
0 & \tx{otherwise} \eec
\]}
where $a = \sum_{i = 0}^k a_i$ and $n = \sum_{i=0}^k x_i$ with $x_0 = \ga$ and $a + (n-1) c > 0$.  This is the well-known inverse
P\'{o}lya-Eggenberger distribution.  A proof of this distribution is given as follows.

To make the sampling experiment well-defined,  assume that for $i = 0, 1, \cd, k$, there are infinitely many balls of color $\bb{C}_i$.  Label
these balls as $B_{i, 1}, B_{i, 2}, \cd$ so that $B_{i, j}, \; j = 1, \cd, a_i$ will be initially put in the urn and that balls $B_{i, j}, \; j
= a_i + 1, \; a_i + 2, \cd$ are used as additional balls, which are used in an order consistent with their indexes, that is, the ball with index
$j$ must have been added to the urn if the ball with index $j + 1$ is to be added to the urn.  Let the set of balls be denoted by $\{ B_{i, j}
\}$.

Consider the making of a general permutation (without color restriction), denoted by $P_1, P_2, \cd, P_n$, of balls with $n$ trials, where for
$\ell = 1, \cd, n$, the ball $P_\ell$ is drawn in the $\ell$-th trials.  Clearly, $P_\ell$ is picked from a subset of $\{ B_{i, j} \}$ under
certain constraints.   Note that the same $B_{i, j}$ can appear in different position of $P_1 P_2 \cd P_n$.  At the $1$st trial, there are $a$
balls and thus there are $a$ choices. At the $\ell$-th trial, where $1 \le \ell \leq n$, there are $a + (\ell - 1) c$ balls and thus there are
$a + (\ell - 1) c$ choices. At the $n$-th trial, there are $a + (n-1) c$ balls and accordingly there are $a + (n-1) c$ choices. Multiplying the
numbers of these choices gives the total number of all possible permutations created by these $n$ trials, which is $\prod_{\ell = 1}^n [a +
(\ell - 1) c ]$.  For the experiment to be feasible, $\prod_{\ell = 1}^n [ a + (\ell - 1) c ]$ must be positive.

Now consider the number of permutations $P_1, P_2, \cd, P_n$ corresponding to $x_i$ balls of color $\bb{C}_i$ for $i = 0, 1, \cd, k$. Such
permutations can be made as follows.

First, note that the draw at the $n$-th trial must be of color $\bb{C}_0$. For $i = 0, 1, \cd, k$, we choose $x_i$ numbers from $\{1, 2, \cd, n
- 1\}$ and denote the set of the numbers as $S_i$. Thus, we have $k$ sets $S_0, S_1, \cd, S_k$. There are $\f{(n-1)!}{(x_0 -1)! x_1 !, \cd, x_k
!}$ ways to create such sets.

Second, for $i = 0, 1, \cd, k$, arrange the numbers in $S_i$ as $r_{\ell, i}, \; \ell = 1, \cd, x_i$.  In the $r_{\ell, i}$-th trials, there are
$a_i + (\ell - 1) c$ balls of color $\bb{C}_i$.  Hence, for $i = 0, 1, \cd, k$, there are $\prod_{\ell = 1}^{x_i} [ a_i + (\ell - 1) c ]$ ways
to make a permutation having $x_i$ balls if the pattern of balls of other colors are fixed.  It follows that there are $\f{(n-1)!}{(x_0 -1)! x_1
!, \cd, x_k !} \prod_{i=0}^k \li \{ \prod_{\ell = 1}^{x_i} [a_i + (\ell - 1) c ] \ri \}$ permutations having $x_i$ balls of color $\bb{C}_i$ for
$i = 0, 1, \cd, k$. For the experiment to be feasible, $\prod_{i=0}^k \li \{ \prod_{\ell = 1}^{x_i} [a_i + (\ell - 1) c ] \ri \}$ must be
positive.

Third, since each permutation is equally likely, the probability of getting $x_i$ balls of color $\bb{C}_i$ respectively is equal to the ratio
of the number of permutations, having $x_i$ balls of color $\bb{C}_i$ for $i = 0, 1, \cd, k$,  to the total number of permutations, that is,
\[
\f{1}{ \prod_{\ell = 1}^n [a + (\ell - 1) c ] } \f{(n-1)!}{(x_0 -1)! x_1 !, \cd, x_k !} \prod_{i=0}^k \li \{ \prod_{\ell = 1}^{x_i} [a_i + (\ell
- 1) c ] \ri \},
\]
the probability mass function is thus justified.

\end{document}